\documentclass[12pt]{amsart}
\usepackage{mathtools,latexsym,amssymb,amsthm,amsmath,enumerate,epsfig,setspace,bbding,color,graphicx,caption,epstopdf,thmtools,scrextend}
\usepackage[all]{xy}
\usepackage[retainorgcmds]{IEEEtrantools}
\usepackage{hyperref}
\usepackage[margin=10pt,font=small,labelfont=bf]{caption}
\usepackage{tikz}
\usetikzlibrary{matrix,arrows.meta,decorations.pathmorphing,positioning}

\definecolor{darkred}{rgb}{0.5,0,0}
\definecolor{darkgreen}{rgb}{0,0.5,0}
\definecolor{darkblue}{rgb}{0,0,0.5}
\definecolor{darkorange}{rgb}{0.3,0.6,0.2}
\definecolor{darkyellow}{rgb}{0.75,0.75,0.2}

\renewcommand{\phi}{\varphi}
\numberwithin{equation}{section}
\theoremstyle{plain}

\declaretheorem[name={Theorem},style=plain,numberwithin=section]{thm}

\declaretheorem[name={Corollary},style=plain,sibling=thm]{cor}
\declaretheorem[name={Proposition},style=plain,sibling=thm]{prop}

\declaretheorem[name={Definition},style=definition,qed=$\diamondsuit$,sibling=thm]{defi}

\declaretheorem[name={Remark},style=definition,sibling=thm]{rmk}
\declaretheorem[numbered=no,name={Theorem},style=plain]{thm*}
\declaretheorem[numbered=no,name={Acknowledgements},style=definition]{ack}
\declaretheorem[numbered=no,name={Main Theorem},style=plain]{mainthm}

\newcommand{\U}{\text{$\mathcal{U}$}}

\newcommand{\R}{\mathbb{R}}
\newcommand{\C}{\mathbb{C}}
\newcommand{\T}{\mathbb{T}}
\newcommand{\Z}{\mathbb{Z}}

\newcommand{\X}{\mathfrak{X}}
\newcommand{\g}{\mathfrak{g}}

\newcommand{\id}{\mathrm{id}}

\newcommand{\im}{\mathrm{im}}

\newcommand{\pr}{\mathrm{pr}}

\newcommand{\G}{\mathcal{G}}

\newcommand{\F}{\mathcal{F}}
\renewcommand{\U}{\mathcal{U}}
\renewcommand{\H}{\mathcal{H}}
\newcommand{\s}{\mathbf{s}}

\renewcommand{\t}{\mathbf{t}}
\renewcommand\phi{\varphi}
\renewcommand{\P}{\mathbb{P}}
\newcommand\twiddle[1]{{\widetilde{#1}}}
\newcommand{\comment}[1]{}

\newcounter{prof}
\newenvironment{prof}[1][Proof.]
  {\stepcounter{prof}\begin{proof}[#1]}
  {\end{proof}}

\makeatletter
\@addtoreset{claim}{prof}
\makeatother

\title{Poisson manifolds of strong compact type over 2-tori}
\date{September 1, 2023}
\author{Luka Zwaan}
\thanks{This work was partially supported by NSF grant DMS-2003223.}
\begin{document}

\begin{abstract}
We construct a new class of examples of Poisson manifolds of strong compact type. In particular, we show that all strongly integral affine circles and two-dimensional tori appear as the leaf space of a Poisson manifold of strong compact type. 
\end{abstract}

\maketitle

\section{Introduction}

Like symplectic geometry, Poisson geometry started from the mathematical formalisation of classical mechanics. Roughly speaking, a Poisson manifold is a smooth manifold equipped with a \emph{Poisson bracket} on its space of smooth functions, which allows one to formulate Hamiltonian dynamics. Examples of Poisson manifolds include symplectic manifolds and duals of Lie algebras, an early glimpse into the deep connection with symplectic geometry and Lie theory. Unlike symplectic manifolds, Poisson manifolds are very flexible in nature. For instance, every manifold admits a Poisson structure and there is no local classification of Poisson structures. For this reason it is common to restrict one's attention to specific classes of Poisson manifolds, where one can formulate deep results about their geometry. In this paper we are concerned with \emph{Poisson manifolds of compact type} (PMCTs). PMCTs are the ``compact objects'' in Poisson geometry. They were first introduced in \cite{rigidity} and their role in the theory is analogous to the one played by compact Lie algebras in Lie theory. Just as there is the special class of compact \emph{semisimple} Lie algebras among compact Lie algebras, there is an important distinguished class among PMCTs, namely that of Poisson manifolds of \emph{strong} compact type (PMSCTs). A simple class of examples of PMSCTs is given by compact symplectic manifolds with finite fundamental group, but it is difficult to construct examples that are not symplectic. The first such example was given in \cite{2013arXiv1312.7267M}, building on work of \cite{kotschick2006free}. There a regular PMSCT is constructed whose symplectic leaves are all diffeomorphic to a K3 surface and whose leaf space is diffeomorphic to a circle. One can form new PMSCTs by taking products of the aforementioned examples, but apart from these no other examples are known. In this paper we use the construction of \cite{2013arXiv1312.7267M} to obtain new examples of PMSCTs. It is known that the leaf space of a PMSCT must be a compact integral affine orbifold and in the example of \cite{2013arXiv1312.7267M} this is the ``standard'' integral affine structure on the circle. In this work we show that \emph{all} strongly integral affine circles and two-dimensional tori can appear as the leaf space of a PMSCT.

In order to explain our main result, recall that a Poisson structure on a manifold $M$ is a Lie bracket on $C^{\infty}(M)$ which is a derivation in each entry. Equivalently, a Poisson structure is a bivector $\pi\in\X^2(M)$ satisfying $[\pi,\pi]=0$. This is the definition we work with in this paper. Every Poisson manifold has a partition into symplectic manifolds. This \emph{symplectic foliation} can be viewed as a singular foliation integrating the (singular) distribution $\pi^{\#}(T^*M)\subset TM$. If $\pi$ has constant rank, this is actually a regular foliation. In this case the Poisson manifold is called \emph{regular}.

The ``global'' objects in Poisson geometry are the so-called \emph{symplectic groupoids}. A symplectic groupoid is a Lie groupoid $\G\rightrightarrows M$ carrying a multiplicative symplectic form $\Omega\in\Omega^2(\G)$. A Poisson manifold $(M,\pi)$ is called \emph{integrable} if there exists some symplectic groupoid $(\G\rightrightarrows M,\Omega)$ for which the target map $\t:(\G,\Omega)\to (M,\pi)$ is a Poisson map (see \cite{lectures}). PMCTs are defined as those Poisson manifolds that are integrated by a source connected, Hausdorff symplectic groupoid having a certain compactness property. Contrary to the case of Lie groups and Lie algebras, there are multiple notions of compactness for Lie groupoids, namely a Lie groupoid $\G\rightrightarrows M$ is called
\begin{itemize}
\item \emph{proper} if the anchor map $(\s,\t):\G\to M\times M$ is proper;
\item \emph{source proper}, or $\s$-proper, if the source map is proper;
\item \emph{compact} if the space of arrows $\G$ is compact.
\end{itemize}
Accordingly, we say that $(M,\pi)$ is of \emph{proper/source proper/compact type} if it admits a source connected, Hausdorff symplectic groupoid of proper/source proper/compact type, respectively.

The types just defined depend on the choice of integration of $(M,\pi)$. However, just like for Lie groups, there is a unique ``largest'' integration, namely the one with 1-connected source fibers. This is often called the \emph{Weinstein groupoid}. We say that an integrable Poisson manifold has \emph{strong} proper/source proper/compact type if its Weinstein groupoid is Hausdorff and has the corresponding type. As mentioned above, we will focus here on Poisson manifolds of strong compact type.

Unlike general Poisson manifolds, PMCTs have a rich geometry transverse to their associated symplectic foliation. For example, the leaf space of a regular PMCT inherits the structure of an integral affine orbifold. Roughly speaking this means that the leaf space has an orbifold atlas where the transitions are integral affine maps. The precise statement can be found in \cite{pmct1,pmct2}, where many other properties of PMCTs are discussed.

As mentioned above, the first example of a PMSCT that is not symplectic was given by Martinez-Torres in \cite{2013arXiv1312.7267M}. The construction there is inspired by the work of Kotschick \cite{kotschick2006free}, where non-trivial results on the geometry of K3 surfaces are used to construct a free symplectic circle action with contractible orbits. The orbit space of such an action is a PMSCT with smooth leaf space a circle endowed with its standard integral affine structure (that is, the one it inherits as a quotient of $\R$ by $\Z$ acting by translations). In general, it is not known whether any compact integral affine orbifold can appear as the leaf space of a PMSCT. On the one hand constructing strong PMCTs is a difficult problem on its own, and on the other not much is known about the classification of compact integral affine manifolds in dimension greater than two. The integral affine structures on a circle are easily classified, and the classification of integral affine structures on compact 2-dimensional manifolds was obtained in \cite{MISHACHEV1996301,Kleinintaff}. The main result of this paper is the following.

\begin{mainthm}
Any strongly integral affine circle or two-dimensional torus can be realised as the leaf space of a PMSCT.
\end{mainthm}

Here by a \emph{strongly integral affine structure} we mean an integral affine structure with integral translational part (see \cite[Remark 5.10]{strongintaff} and Remark \ref{rmk:strongly}).

Our strategy to prove this result is as follows. Using the geometry of K3 surfaces one constructs a universal family of marked K\"ahler K3 surfaces (see Section \ref{sec:k3}) to which one can apply a general method from \cite{pmct2} to obtain PMSCTs. Using this construction together with the classification of integral affine 2-tori from \cite{MISHACHEV1996301}, one obtains examples of PMSCTs for all isomorphism classes of strongly integral affine 2-tori.

This paper is organised as follows. In Section \ref{sec:gen}, we provide some background on PMCTs and we recall the general method of constructing regular PMSCTs from \cite{pmct2}. In Section \ref{sec:k3} we recall the relevant results on K3 surfaces that are needed for our construction. The resulting examples of PMSCTs have symplectic foliation a fibration over $S^1$ or $\T^2$ with typical fibre the smooth manifold underlying a K3 surface. The symplectic structures on the fibres vary in a controlled fashion which ensures that the Weinstein groupoid is a compact symplectic groupoid.
Finally, Section \ref{sec:ex} is dedicated to the actual constructions, which includes some lengthy computations. We treat the circle case first and this includes the original example from \cite{2013arXiv1312.7267M}. Lastly, we construct the PMSCTs with leaf space the strongly integral affine 2-tori.

\begin{ack}
I would like to thank Marius Crainic for his help as supervisor of my master thesis, when I first worked on this topic. I would also like to thank David Mart\'{i}nez Torres for a fruitful discussion regarding his paper that this work is based on. Finally, I would like to thank my PhD advisor Rui Loja Fernandes for his feedback during the writing of this paper.
\end{ack}

\section{Background \& general construction of PMSCTs}\label{sec:gen}

The construction we give below is based on two results on PMCTs:
\begin{enumerate}[(a)]
\item the leaf space carries an integral affine orbifold structure (see \cite[Section 3]{pmct2}), and
\item the \emph{linear variation theorem} (see \cite[Section 4-5]{pmct2}).
\end{enumerate}

We briefly recall these results before giving the general construction. Here we only need to consider the case of $1$-connected leaves. In this case the leaf space is smooth, since this assumption implies that the monodromy groupoid of the symplectic foliation is proper and has trivial isotropy groups. Then both (a) and (b) above simplify significantly.

\subsection{The integral affine structure on the leaf space}\label{sec:intaff}

Recall that an integral affine structure on a manifold $B$ is given by an atlas whose transition functions are integral affine maps. Equivalently, it is specified by a lattice $\Lambda\subset T^*B$ locally spanned by closed 1-forms. 

Consider a regular, $\s$-connected, proper symplectic groupoid $(\G,\Omega)\rightrightarrows(M,\pi)$. As mentioned above, we assume that the associated symplectic foliation $\F_{\pi}$ has 1-connected leaves so that the leaf space $B$ is a smooth manifold. We obtain a lattice $\twiddle{\Lambda}\subset\nu^*(\F_{\pi})$ as follows:
\begin{enumerate}
\item for each $x\in M$, the kernel of the exponential map $\g_x\to\G_x$ gives a lattice in $\g_x$ and
\item the isomorphism $\g_x\cong\nu_x^*(\F_{\pi})$ induced by $\Omega$ allows us to transport it to the conormal space.
\end{enumerate}
This lattice descends to an integral affine structure $\Lambda\subset T^*B$ on $B$.

\subsection{The linear variation theorem}\label{sec:linvar}

We asssume now in addition that $(\G,\Omega)\rightrightarrows(M,\pi)$ is source proper. Denoting the symplectic leaf corresponding to $b\in B$ by $(S_b,\omega_b)$, we form the vector bundle
\[ \H^2:=\bigsqcup_{b\in B} H^2(S_b,\R)\to B
\]
and the lattice
\[ \H^2_{\Z}:=\bigsqcup_{b\in B} \im\big(H^2(S_b,\Z)\to H^2(S_b,\R)\big)
\]
inside it. Associated to this we have the \emph{Gauss-Manin connection} $\nabla$ on $\H^2$, uniquely determined by requiring the sections of $\H^2_{\Z}$ to be parallel. Note that $\pi$ gives us a section $\varpi\in\Gamma(\H^2)$, $b\mapsto[\omega_b]$.

The Gauss-Manin allows us to study the variation of $\varpi$: parallel transport makes $\H^2$ into a $\Pi_1(B)$-representation and we define the \emph{variation map} $\mathrm{var}_{\varpi}:\Pi_1(B)\to\H^2$ to be
\[ [\gamma]\mapsto \gamma_*(\varpi_{\gamma(0)}) \in \H^2_{\gamma(1)}.
\] 
On the other hand, we also have the \emph{linear variation map} $\mathrm{var}_{\varpi}^{\mathrm{lin}}:TB\to\H^2$ given by
\[ v\mapsto\nabla_v\varpi
\]
and the \emph{affine variation map} $\mathrm{var}_{\varpi}^{\mathrm{aff}}:=\varpi+\mathrm{var}_{\varpi}^{\mathrm{lin}}$.

The linear variation theorem relates the variation and affine variation maps by means of the \emph{developing map} associated to the integral affine manifold $(B,\Lambda)$.  Associated to the lattice $\Lambda^*\subset TB$ we have a canonical flat connection on $TB$ (not to be confused with $\nabla$ above). This makes $TB$ into a $TB$-representation, and since the connection is torsion-free the identity map $TB\to TB$ is an algebroid cocycle. The developing map is defined to be the groupoid cocycle $\mathrm{dev}:\Pi_1(B)\to TB$ integrating it.

\begin{rmk}
One can show that after fixing $b\in B$ and a basis of $\Lambda_b$ this boils down to the classical notion of developing map defined on the universal covering space (see \cite[Section 4.2]{pmct2}):
\[ \mathrm{dev}_b:\widetilde{B}\to T_b B\simeq\R^q.\] 
\end{rmk}

We can now state the linear variation theorem as follows.

\begin{thm}[{{\cite[Theorem 4.4.2]{pmct2}}}]\label{thm:linvar}
One has a commutative diagram
\begin{center}
\begin{tikzpicture}
\matrix(m)[matrix of math nodes,
row sep=3em, column sep=4em,
text height=1.5ex, text depth=0.25ex]
{
\Pi_1(B) & & \H^2 \\
 & TB & \\
};
\begin{scope}[every node/.style={midway,auto,font=\small}]
\draw[->] (m-1-1) -- node {$\mathrm{var}_{\varpi}$} (m-1-3);
\draw[->] (m-1-1) -- node[swap] {$\mathrm{dev}$} (m-2-2);
\draw[->] (m-2-2) -- node[swap] {$\mathrm{var}_{\varpi}^{\mathrm{aff}}$} (m-1-3);
\end{scope}
\end{tikzpicture}
\end{center}
\end{thm}

This rather abstract formulation can locally be made explicit. Let $b_0\in B$ and choose an integral affine chart $(U,\varphi)$ centered at $b_0$ such that $\varphi(U)$ is convex and such that $M\to B$ trivialises over $U$. This induces a trivialisation $\Phi:\H^2|_U\cong U\times H^2(S_{b_0},\R)$. The chart induces an identification $T_{b_0}B\cong\R^q$ and allows us to consider ``straight line'' paths from $b\in U$ to $b_0$. Restricting to such paths the above diagram becomes
\begin{center}
\begin{tikzpicture}
\matrix(m)[matrix of math nodes,
row sep=3em, column sep=4em,
text height=1.5ex, text depth=0.25ex]
{
U & & H^2(S_{b_0},\R) \\
 & \R^q & \\
};
\begin{scope}[every node/.style={midway,auto,font=\small}]
\draw[->] (m-1-1) -- node {$b\mapsto\Phi([\omega_b])$} (m-1-3);
\draw[->] (m-1-1) -- node[swap] {$\varphi$} (m-2-2);
\draw[->] (m-2-2) -- node[swap] {$v\mapsto[\omega_{b_0}]+\sum_iv_ic_i$} (m-1-3);
\end{scope}
\end{tikzpicture}
\end{center}
where $c_i\in H^2(S_{b_0},\Z)$ are the Chern classes of the torus bundle $\s^{-1}(x)\to S_{b_0}$, where $x\in S_{b_0}$ (see \cite[Corollary 4.4.4]{pmct2}). This local formulation is reminiscent of the linear variation theorem from \cite{dh}. In other words, Theorem \ref{thm:linvar} can be viewed as a global formulation and generalisation of the classical Duistermaat-Heckman theorem.

\subsection{The construction} 
The construction we describe in this section yields a PMSCT with 1-connected symplectic leaves, whose leaf space is a complete integral affine manifold. This means that the leaf space is a quotient of $\R^q$ by a free and proper action of a discrete group of integral affine transformations. Note that if the Markus conjecture holds true, then in fact every compact integral affine manifold is of this type (see \cite[Section 8.6]{Goldman22}). This allows us to give an explicit formulation of the linear variation, similar to the discussion following Theorem \ref{thm:linvar}.
The setup is as follows.

Let $E\to\R^q$ be a fibre bundle with typical fibre $S$, a compact 1-connected manifold, and assume that $E$ admits a Poisson structure $\pi_E$ whose symplectic leaves are precisely the fibres of this bundle. As in Section \ref{sec:linvar} we have (i) the vector bundle $\H^2\to\R^q$ whose fibers are the degree two cohomology groups of the symplectic leaves, (ii) the lattice $\H^2_{\Z}\subset\H^2$ of integral cohomology, (iii) the associated Gauss-Manin connection $\nabla$ and (iv) the section $\varpi\in\Gamma(\H^2)$ induced by $\pi_E$.

Next, let $\Gamma\subset\mathrm{Aff}_{\Z}(\R^q)=\{x\mapsto Ax+v\mid A\in \mathrm{GL}(q,\Z), v\in \R^q\}$ be a discrete group of integral affine transformations acting freely and properly on $\R^q$, and assume that there is a Poisson action of $\Gamma$ on $(E,\pi_E)$ making the projection $E\to\R^q$ equivariant. Then setting $M:=E/\Gamma$ and $B:=\R^q/\Gamma$, we get a (smooth) fibre bundle $p:M\to B$, again with typical fibre $S$, and a Poisson structure $\pi$ on $M$ whose leaves are the fibres of $p$. In other words, $(M,\pi)$ is a regular Poisson manifold with leaf space $B$. Note also that $B$, being a quotient $\R^q/\Gamma$, naturally inherits an integral affine structure.

We can now state the general method of constructing PMSCTs. It is a reformulation of \cite[Proposition 4.4.6]{pmct2}. 

\begin{prop}\label{prop:genconstr}
Let $(M=E/\Gamma,\pi)$ be constructed as above. Assume that there exists a $\nabla$-flat section $s\in\Gamma(\H^2)$ and linearly independent sections $c_1,\ldots,c_q\in\Gamma(\H^2_{\Z})$ such that
\begin{equation}\label{eq:linvar}
\varpi=s+\sum_{i=1}^q \pr^i\cdot c_i,
\end{equation}
where $\pr^i:\R^q\to\R$ denotes projection onto the $i$-th coordinate. Then $(M,\pi)$ is of strong $\s$-proper type and the induced integral affine structure on $B$ agrees with the one coming from the quotient $\R^q/\Gamma$. In particular, if $B$ is compact then $(M,\pi)$ is a PMSCT.
\end{prop}

\begin{prof}
Pulling back the integral affine structure on $B$ along $p:M\to B$ yields a transverse integral affine structure on the symplectic foliation $\F_{\pi}$, \emph{i.e.}~a lattice in its conormal bundle. We denote this lattice by $\twiddle{\Lambda}\subset\nu^*(\F_{\pi})$. The main point is that for all $x\in M$, the monodromy group $N_x(M,\pi)$ is equal to the lattice $\twiddle{\Lambda}_x$. In fact, using the description of the monodromy groups for regular Poisson manifolds as the variation of symplectic areas (see \cite[Section 6]{intpoiss}) this follows directly from equation (\ref{eq:linvar}). The integrability criteria for Poisson manifolds then imply that $(M,\pi)$ is integrable. Furthermore, since $S$ has trivial fundamental group, the isotropy groups of the Weinstein groupoid $\Sigma(M,\pi)$ fit into the exact sequence
\[ \cdots\to\pi_2(S,x)\xrightarrow{\partial_x}\nu^*_x(\F_{\pi})\to\Sigma_x(M,\pi)\to0,
\]
where $\partial_x$ is the monodromy map at $x$. Therefore, from our previous discussion, it follows that $\Sigma_x(M,\pi)\simeq\nu^*_x(\F_{\pi})/\twiddle{\Lambda}_x$, \emph{i.e.}~that the isotropy group at $x$ is compact. Since this holds for all $x\in M$ and since $S$ is also compact, this shows that the Weinstein groupoid is $\s$-proper.

Finally, since $\twiddle{\Lambda}\subset\nu^*(\F_{\pi})$ is closed, Hausdorffness of the Weinstein groupoid follows from \cite[Theorem 1.1]{AlcaldeCuesta1994IntegrationSD}.
\end{prof}

\section{Background on K3 surfaces and the Poisson structure on the universal family}\label{sec:k3}

We start by listing some definitions and results concerning K3 surfaces, after which we describe the moduli spaces and universal families for K3 surfaces. These results can be found in \cite{barth1984compact}. Finally, following \cite{2013arXiv1312.7267M}, we use the Calabi-Yau theorem to turn the universal family into a Poisson manifold and the strong Torelli theorem to establish a Poisson action on it, setting us up to apply our construction.

\begin{defi}
A K3 surface is a compact, 1-connected complex surface with trivial canonical bundle.
\end{defi}

Every K3 surface is K\"{a}hler (see \cite{Siu1983}). All K3 surfaces have the same underlying smooth manifold $S$ (see \cite[Corollary VIII.8.6]{barth1984compact}); this will be the model fibre used in Proposition \ref{prop:genconstr}. The intersection form on $H^2(S,\Z)$ turns it into a lattice and this lattice is isomorphic to the aptly named \emph{K3 lattice}, which we denote by $(L,(\cdot,\cdot))$. It is the unique even, unimodular lattice of signature $(3,19)$ (see \cite[Proposition VIII.3.2 (ii)]{barth1984compact}). Explicitly, we have $L=U^{\oplus3}\oplus(-E_8)^{\oplus2}$, where $U=\Z^{\oplus2}$ with form given by $\begin{pmatrix}
0&1\\1&0
\end{pmatrix}$ and $E_8=\Z^{\oplus8}$ with form given by the Cartan matrix of $E_8$; it is important for us that this form is positive definite. We also set $L_{\R}:=L\otimes\R$ and $L_{\C}:=L\otimes\C$; note that these are models for the real and complex cohomology, respectively.

\subsection{The Torelli theorem}

\begin{defi}
Let $X,X'$ be K3 surfaces. A $\Z$-module isomorphism $H^2(X',\Z)\to H^2(X,\Z)$ is a \emph{Hodge isometry} if
\begin{enumerate}[(i)]
\item it preserves the intersection form;
\item its $\C$-linear extension preserves the Hodge decomposition.
\end{enumerate}
A Hodge isometry is called \emph{effective} if its $\R$-linear extension maps some K\"{a}hler class of $X'$ to one of $X$.
\end{defi} 

Effectiveness of a Hodge isometry is equivalent to requiring it to map the K\"{a}hler cone of $X'$ to that of $X$ (see \cite[Proposition VIII.3.10]{barth1984compact}).

\begin{thm}[{{Torelli \cite[Corollary VIII.11.4]{barth1984compact}}}]
\label{thm:torelli}
Let $X,X'$ be K3 surfaces. Then for any effective Hodge isometry $\varphi:H^2(X',\Z)\to H^2(X,\Z)$ there exists a unique biholomorphism $f:X\to X'$ such that $f^*=\varphi$.
\end{thm}

This result is ultimately used to obtain the action in Proposition \ref{prop:genconstr}.

\subsection{Moduli spaces and universal families}

There are two moduli spaces and corresponding families for K3 surfaces: one takes into account the K\"{a}hler structure and the other only considers the complex structure. We start now with the latter. 

\begin{defi}
A \emph{marked K3 surface} is a pair $(X,\varphi)$ consisting of a K3 surface $X$ and a \emph{marking} $\varphi$, \emph{i.e.} an isometry $\varphi:H^2(X,\Z)\to L$. Two marked K3 surfaces are \emph{equivalent} if there is a bihomolorphism between them intertwining the markings. The \emph{moduli space of marked K3 surfaces} is the set of equivalence classes:
\[ M_1:=\{(X,\varphi)\}/\sim. \qedhere
\]
\end{defi}

It follows immediately from the definition that any K3 surface admits, up to scalar multiplication, a unique nowhere vanishing holomorphic $2$-form. In fact, one can show that, again up to scalar multiplication, there is a bijection between complex structures on $S$ and closed, complex 2-forms $\sigma\in\Omega^2(S,\C)$ satisfying $\sigma\wedge\sigma=0$ and $\sigma\wedge\bar{\sigma}>0$. This motivates the following definitions. We will use the same letter to denote a marking $\varphi:H^2(X,\Z)\to L$  and the induced maps $\varphi:H^2(X,\R)\to L_{\R}$ and $\varphi:H^2(X,\C)\to L_{\C}$.

\begin{defi}
The \emph{period domain} is given by
\[ \Omega:=\{[\sigma]\in\P(L_{\C})\mid (\sigma,\sigma)=0, (\sigma,\bar{\sigma})>0\}.
\]
We define the \emph{period map} $\tau_1:M_1\to\Omega$ by
\[ [(X,\varphi)]\mapsto[\varphi(\sigma_X)],
\]
where $\sigma_X$ is a nowhere vanishing holomorphic $2$-form on $X$.
\end{defi}

\begin{thm}[{{\cite[Theorem VIII.12.1]{barth1984compact}}}]
The moduli space $M_1$ admits the structure of a $20$-dimensional complex manifold such that the period map $\tau_1:M_1\to\Omega$ becomes a surjective local biholomorphism. Furthermore, there exists a universal family $\U\to M_1$ of marked K3 surfaces. 
\end{thm}

\begin{rmk} Recall that a family is universal if any other family is locally the pullback of it by a unique map (see \cite[Section I.10]{barth1984compact}). The fibre of the universal family $\U\to M_1$ over any $t\in M_1$ is a marked K3 surface $(X_t,\varphi_t)$ such that $[(X_t,\varphi_t)]=t$. Furthermore, these markings vary smoothly in the sense that they induce local trivialisations of the bundle $\cup_{t\in M_1} H^2(X_t,\R)$.
\end{rmk}

There are still some inconveniences present here. It can be shown that $M_1$ is not Hausdorff, and that the period map $\tau_1$ is not injective (see \cite[Remark VIII.12.2]{barth1984compact}). These problems disappear when taking into account the K\"{a}hler structure. 

\begin{defi}
We define $M_2$ to be the subset of the bundle
\[ \bigsqcup_{t\in M_1} H^2(X_t,\C)
\]
consisting of all K\"{a}hler classes. 
\end{defi}

It can be shown that $M_2$ is a real-analytic manifold of dimension $60$ (see \cite[Lemma VIII.9.3]{barth1984compact} and its proof). One should think of a point in $M_2$ as an equivalence class of marked K3 surfaces together with a specified K\"{a}hler class. Note that there is a projection map $\pr:M_2\to M_1$.

Inspired by some analysis of the K\"{a}hler cone of K3 surfaces (see \cite[Section VIII.3 and VIII.9]{barth1984compact}) one makes the following definitions.

\begin{defi} Set
\[ K\Omega:=\{(k,[\sigma])\in L_{\R}\times\Omega\mid (k,k)>0,(k,\sigma)=0\}.
\]
The \emph{refined period domain} is then given by
\[ K\Omega^0:=\{(k,[\sigma])\in K\Omega\mid (k,d)\neq0 \text{ for all } d\in L \text{ such that } (d,d)=-2 \text{ and } (d,\sigma)=0\}.
\]
The \emph{refined period map} $\tau_2:M_2\to K\Omega^0$ is defined as
\[ (t,k)\mapsto(\varphi_t(k),\tau_1(t)).\qedhere
\]
\end{defi}

\begin{thm}[{{\cite[Theorem VIII.12.3 and VIII.14.1]{barth1984compact}}}]
The refined period map is a diffeomorphism.
\end{thm}

We set $K\U:=(\pr\circ\tau_2^{-1})^*\U$. This is a real-analytic family (\emph{i.e.}~fibre bundle) over $K\Omega^0$ with extra data attached: the fibre over $(k,[\sigma])$ is a triple $(X,\varphi,\omega)$ consisting of a K3 surface $X$, a marking $\varphi:H^2(X,\Z)\to L$ and a K\"{a}hler class $\omega\in H^2(X,\R)$ such that $\varphi(\omega)=k$. These markings vary smoothly in the same sense as before, and hence so do the K\"{a}hler classes.

The family $K\U\to K\Omega^0$ is universal for real-analytic ``marked K\"{a}hler K3 families'', \emph{i.e.}~real-analytic families of K3 surfaces equipped with smoothly varying markings and K\"{a}hler classes.

\subsection{The Poisson structure}

Recall the following special version of the Calabi-Yau theorem (see e.g. \cite[Theorem I.15.1]{barth1984compact}).

\begin{thm}
Let $X$ be a compact complex manifold with vanishing first Chern class. Then for any K\"{a}hler class $\omega\in H^2(X,\R)$ there exists a unique Ricci flat K\"{a}hler metric whose K\"{a}hler form belongs to $\omega$.
\end{thm}

This theorem applies in particular to K3 surfaces, and thus we can use it to endow the fibres of $K\U\to K\Omega^0$ with smoothly varying K\"{a}hler forms, turning it into a Poisson manifold (see also \cite[Secion 2.1.3]{2013arXiv1312.7267M}).

\begin{cor}\label{cor:leafspace}
The family $K\U$ admits a regular Poisson structure $\pi_{K\U}$ whose symplectic leaves are the fibres of $K\U\to K\Omega^0$. Moreover the symplectic form on the fibre $X$ over $(k,[\sigma])$ with marking $\varphi$ is the K\"{a}hler form associated to the unique Ricci flat K\"{a}hler metric on $X$ with K\"{a}hler class $\varphi^{-1}(k)$.
\end{cor}

\subsection{The action}

We will construct an action on $K\U$ by the group $O(L)$ of isometries of the K3 lattice. Note that there is an obvious induced action of $O(L)$ on $K\Omega^0$.

\begin{prop}
\label{prop:action}
There is a Poisson action of $O(L)$ on $(K\U,\pi_{K\U})$ with respect to which the projection $K\U\to K\Omega^0$ is equivariant.
\end{prop}
\begin{prof}
Fix $\gamma\in O(L)$ and $p\in K\Omega^0$. Using the notation from above, denote the triple over $p$ by $(X_p,\varphi_p,\omega_p)$ and similarly for $\gamma(p)$. It is easy to see that
\[ \varphi_p^{-1}\circ\gamma^{-1}\circ\varphi_{\gamma(p)}:H^2(X_{\gamma(p)},\Z)\to H^2(X_p,\Z)
\]
is an effective Hodge isometry, so that by Theorem \ref{thm:torelli} we obtain a biholomorphism $f_{\gamma}^p:X_p\to X_{\gamma(p)}$. The universality of the family then gives neighbourhoods $U$ and $V$ of $p$ and of $\gamma(p)$ respectively and an isomorphism $(\Psi,\psi):K\U|_U\to K\U|_V$ extending $f_{\gamma}^p$: through the biholomorphism $f_{\gamma}^p$, $K\U$ becomes a deformation of $X_p$ at two basepoints, $p$ and $\gamma(p)$. Since $K\U$ is universal, these two deformations are locally isomorphic. Writing $\Psi_q:X_q\to X_{\psi(q)}$ for the fiberwise maps, it then follows that for all $q\in U$ we have that
\[ \Psi_q^*=\varphi_q^{-1}\circ\gamma^{-1}\circ\varphi_{\psi(q)}:H^2(X_{\psi(q)},\Z)\to H^2(X_q,\Z).
\]
This implies first of all that $\psi=\gamma|_U$, from which it follows that $\Psi_q=f_{\gamma}^q$, since biholomorphisms of K3 surfaces are uniquely determined by their induced maps on degree 2 integral cohomology (see \cite[Proposition VIII.11.3]{barth1984compact}). Thus these fibrewise biholomorphisms $f_{\gamma}^p$, $p\in K\Omega^0$, together form an automorphism $F_{\gamma}:K\U\to K\U$. It is immediate from the above construction that $F_{\id}=\id$, and from the uniqueness part of Theorem \ref{thm:torelli} it follows that $F_{\gamma\circ\gamma'}=F_{\gamma}\circ F_{\gamma'}$ for all $\gamma,\gamma'\in O(L)$, meaning that we have an action of $O(L)$ on $K\U$. This action makes $K\U\to K\Omega^0$ equivariant by construction. Finally, from the uniqueness part of the Calabi-Yau theorem it follows that each $f_{\gamma}^p$ preserves the symplectic forms on the fibres, meaning that the action is by Poisson maps.
\end{prof}

\section{The examples}\label{sec:ex}

From our work in Section \ref{sec:k3} we have a Poisson manifold $(K\U,\pi_{K\U})$ with leaf space $K\Omega^0$ such that:
\begin{enumerate}[(i)]
    \item the cohomology classes of the symplectic forms on the leaves are described in terms of the leaf space $K\Omega^0$ (Corollary \ref{cor:leafspace});
    \item the natural action of $O(L)$ on $K\Omega^0$ lifts to a Poisson action on $(K\U,\pi_{K\U})$ (Proposition \ref{prop:action}).
\end{enumerate}
In order to apply the construction described in Section \ref{sec:gen}, we need to find a suitable embedding $\R^q\hookrightarrow K\Omega^0$ and a suitable subgroup $\Gamma\subset O(L)$. We rephrase Proposition \ref{prop:genconstr} in the current setting in order to make this more precise. For a different version of this result see also \cite[Theorem 1]{2013arXiv1312.7267M}.

\begin{cor}\label{prop:cons}
Assume that we have an embedding $f:\R^q\to K\Omega^0$ and a subgroup $\Gamma\subset O(L)$ such that
\begin{enumerate}[(i)]
\item\label{prop:cons:1} there exist $a\in L_{\R}$ and linearly independent $a_1,\ldots,a_q\in L$ such that the $L_{\R}$-component of $f$ has the form
\[ (x_1,\ldots,x_q)\mapsto a+ \sum_{i=1}^qx_ia_i;
\]
\item\label{prop:cons:2} the action of $\Gamma$ on $K\Omega^0$ preserves the image of $f$;
\item\label{prop:cons:3} the induced action on $\R^q$ is free, proper and by integral affine maps.
\end{enumerate}
Then $M:=f^*K\U/\Gamma$ with the Poisson structure induced from $\pi_{K\U}$ is a Poisson manifold of strong $\s$-proper type with leaf space $B:=\R^q/\Gamma$. If $B$ is compact, $M$ is a PMSCT.
\end{cor}

\begin{rmk}\label{rmk:strongly}
We can now explain why our construction leads to PMSCTs with \emph{strongly} integral affine leaf spaces. 
On the one hand, because of Theorem \ref{thm:linvar}, we are forced to consider embeddings with integral variation, \emph{i.e.}~the $a_i$ must lie in the integral lattice $L$. On the other hand, to apply Theorem \ref{thm:torelli} we need to consider isometries of integral cohomology, \emph{i.e.}~we need to act by elements of $O(L)$. These two technical limitations together only allow for strongly integral affine leaf spaces in the examples.
\end{rmk}

\begin{rmk}
At the level of the symplectic groupoid, one can see that the leaf space being strongly integral affine implies that the restriction of the symplectic form to the identity component of the isotropy (a torus bundle) lies in the integral cohomology. See \cite[Remark 5.10]{strongintaff}.
\end{rmk}

We now recall the classification of strongly integral affine structures for $S^1$ and $\T^2$.

\begin{thm}
    The strongly integral affine circles are, up to isomorphism, the quotients $\R/\Z$ where the $\Z$-action is generated by $x\mapsto x+p$, for a fixed $p\in\Z_{\geq1}$.
\end{thm}
\begin{prof}
    It is easy to see that all integral affine circles are complete. Hence, it suffices to classify, up to conjugation, embeddings $\Z\to\mathrm{Aff}_{\Z}(\R)$ inducing free and proper actions. These are precisely the actions generated by $x\mapsto x+a$ with $a>0$. Restricting to \emph{strongly} integral affine circles yields the result.
\end{prof}

\begin{thm}
    The strongly integral affine 2-tori, up to isomorphism, are quotients $\R^2/\Z^2$, where the $\Z^2$-actions fall into one of the following types:
    \begin{enumerate}[(I)]
    \item\label{torusI} an action generated by $(x,y)\mapsto(x+p,y)$ and $(x,y)\mapsto(x,y+q)$, where $p,q\in\Z_{\geq1}$ and $p|q$;
    \item\label{torusII} an action generated generated by $(x,y)\mapsto(x+p,y)$ and $(x,y)\mapsto(x+ny,y+q)$, where $n,p,q\in\Z_{\geq1}$.
    \end{enumerate}
\end{thm}

\begin{prof}
    The classification of all integral affine structures on 2-tori is given in \cite[Theorem A]{MISHACHEV1996301}. Restricting to strongly integral affine structures and using the Smith normal form for matrices with integer entries to simplify the possibilities from type (\ref{torusI}) yields the above classification.
\end{prof}

\begin{rmk}
The integral affine 2-tori of type (\ref{torusI}) are (isomorphic to) products of integral affine circles. Thus to find examples of PMSCTs with leaf space of this type one can simply take products of PMSCTs with leaf space $S^1$, constructed in Section \ref{example:s1}. This yields Poisson manifolds of dimension 10 whose leaves are products of K3 surfaces. However, the examples we construct in Section \ref{example:t2stan} are six-dimensional Poisson manifolds with K3 surfaces as symplectic leaves and thus result in ``smaller'' examples.
\end{rmk}

\begin{rmk}
Continuing the previous remark, note that by taking products we can also realise \emph{some} higher dimensional integral affine tori as the leaf space of a PMSCT, namely those that are isomorphic to a product of some of the integral affine circles and 2-tori classified above.
\end{rmk}

Before we move on to the examples, we establish some notation. Recall that $L=U^{\oplus3}\oplus(-E_8)^{\oplus2}$. We denote the standard bases of the three copies of $U$ by $\{u,v\}$, $\{x,y\}$ and $\{z,t\}$, so that $(u,v)=(x,y)=(z,t)=1$ with all other combinations yielding zero. Recall also that $-E_8$ is even and negative definite. Finally, let $\{e_1,\ldots,e_8\}$ be a set of real numbers such that the set 
\[ \{1,e_1,\ldots,e_8,e_1^2,e_1e_2,\ldots,e_7^2,e_7e_8,e_8^2\}
\]
consisting of $1,e_1,\ldots,e_8$ and their pairwise products is linearly independent over the integers, or equivalently the rationals. The existence of such a set is guaranteed by \cite{pjm/1103051334}. We then set $e:=(e_1,\ldots,e_8)\in(-E_8)_{\R}$, scaling if necessary such that $|(e,e)|\leq\frac12$, and we set $a:=(0,e), b:=(e,0)\in(-E_8)_{\R}^{\oplus2}\subset L_{\R}$.

Let us outline the strategy for the examples below. In each case, we start by defining $f$ and $\Gamma$. It is fairly straightforward to check items (\ref{prop:cons:2}) and (\ref{prop:cons:3}) from Corollary \ref{prop:cons} and that the image of $f$ is contained in $K\Omega$. It then remains to show that it is actually contained in $K\Omega^0$. This is the more involved part of the computations.

\subsection{The PMSCTs with leaf space the circle}\label{example:s1}
We will construct a PMSCT whose leaf space is a strongly integral affine circle, \emph{i.e.} we want the action of $\Z$ on $\R$ generated by $x\mapsto x+p$ with $p\in\Z_{\geq1}$. The case $p=1$ is the one treated in \cite{2013arXiv1312.7267M} and the computations carried out below for general $p$ are an obvious generalisation of the computations there.

Consider the map $f:\R\to L_{\R}\times\P(L_{\C})$ defined by
\[ s\mapsto (2u+v+sy,[x-su+2y+a+i(z+2t+b)])
\]
and the map $\varphi:L\to L$ defined by $u\mapsto u,v\mapsto v+py$, $x\mapsto x-pu,y\mapsto y$ on the first two copies of $U$ and as the identity on the other summands of $L$. It is easily checked that $\varphi$ is an isometry and that
\[ \varphi\cdot f(s)=f(s+p).
\]
This implies that the image of $f$ is invariant under the action of $\Gamma:=\langle\varphi\rangle$, and also that the induced action on $\R$ is the one we need. 

To show that the image of $f$ is contained in $K\Omega$, let $s\in\R$. Setting $f_1(s)=2u+v+sy$, $f_2(s)=x-su+2y+a$ and $f_3(s)=z+2t+b$, we see that
\begin{align*}
(f_2(s),f_2(s))&=(x-su+2y+a,x-su+2y+a) \\
&=4(x,y)+(a,a) \\
&=4+(e,e) \geq 3\frac12>0, \\
(f_3(s),f_3(s))&=(z+2t+b,z+2t+b) \\
&=4(z,t)+(b,b) \\
&=4+(e,e) \geq 3\frac12>0, \\
(f_2(s),f_3(s))&= (x-su+2y+a,z+2t+b) \\
&=0.
\end{align*}
These computations imply that $[f_2(s)+if_3(s)]\in\Omega$. Since 
\begin{align*}
(f_1(s),f_1(s))&=(2u+v+sy,2u+v+sy)=(2u,v)+(v,2u)=4>0 , \\
(f_1(s),f_2(s))&=(2u+v+sy,x-su+2y+a)=-s(v,u)+s(y,x)=-s+s=0 , \\
(f_1(s),f_3(s))&=(2u+v+sy,z+2t+b)=0 ,
\end{align*}
we see that $f(s)\in K\Omega$. 

It remains to check that $f(s)\in K\Omega^0$ for all $s\in\R$.
\begin{prof}
Assume that we have $d\in L$ such that $(d,d)=-2$ and $(d,f_1(s))=(d,f_2(s))=(d,f_3(s))=0$. We need to find a contradiction. Let us write
\[ d=Au+Bv+Cx+Dy+Ez+Ft+d_1+d_2,
\]
with $A,\ldots,F\in\Z$ and $d_i$ in the $i$-th copy of $-E_8$. Since $E_8$ is even and positive definite, we can write $(d_i,d_i)=-2n_i$, for $n_i\in\Z_{\geq0}$. The above conditions then translate into three equations:
\begin{align}
AB+CD+EF&=n_1+n_2-1,\label{eq:s1:1} \\ 
2B+A+Cs&=0,\label{eq:s1:2} \\ 
D-Bs+2C+(d_2,e)&=0, \label{eq:s1:3} \\ 
F+2E+(d_1,e)&=0. \label{eq:s1:4}
\end{align}

This is where the seemingly strange choice of $e$ comes in. There exist $k_1,\ldots,k_8\in\Z$ such that $(d_1,e)=\sum_ik_ie_i$ and since $\{1,e_1,\ldots,e_8\}$ is linearly independent over the integers by choice of $e$, it follows from (\ref{eq:s1:4}) that we must have $F+2E=k_1=\cdots=k_8=0$. Since the bilinear form on $-E_8$ is nondegenerate, it follows that $d_1=0$ and thus that $n_1=0$.
\smallskip

\textbf{Case $C=0$:} Equation (\ref{eq:s1:2}) yields $2B+A=0$, and (\ref{eq:s1:1}) becomes 
\[ 2B^2+2E^2=1-n_2.
\]
This implies that $B=E=0$ and $n_2=1$. But then $d_2\neq0$ and (\ref{eq:s1:3}) becomes
\[ D+(d_2,e)=0,
\]
which together with $d_2\neq0$ contradicts the ``linear independence'' assumption on $e$.
\smallskip

\textbf{Case $C\neq0$:} From (\ref{eq:s1:2}) we get 
\[ s=-\frac{2B+A}{C},
\]
and substituting this into (\ref{eq:s1:3}) yields
\[ AB+CD=-2C^2-2B^2-(d_2,e).
\]
Combining this with (\ref{eq:s1:1}) gives 
\[ 2B^2+2C^2+2E^2+C(d_2,e)=1-n_2.
\]
From the properties of $e$ we get $Cd_2=0$, implying that $d_2=0$ and thus also that $n_2=0$, so that we are left with
\[2B^2+2C^2+2E^2=1,
\]
which is absurd since $B,C,E\in\Z$.
\end{prof}

\subsection{The PMSCTs with leaf space a torus of type (\ref{torusI})}\label{example:t2stan}

Here we construct a PMSCT with leaf space the torus $\T^2$ with an integral affine structure of type (\ref{torusI}). This means that we want the action of $\Z^2$ on $\R^2$ generated by $(x,y)\mapsto(x+p,y)$ and $(x,y)\mapsto(x,y+q)$, with $p,q\in\Z_{\geq1}$. 

Consider the map $f:\R^2\to L_{\R}\times\P(L_{\C})$ defined by
\[ (s,r)\mapsto (2u+v+sy+rt,[x-su+2y+a+i(z-ru+2t+b)]),
\]
the map $\varphi:L\to L$ as in the previous example and the map $\psi:L\to L$ defined by $u\mapsto u,v\mapsto v+qt$, $x\mapsto x$, $y\mapsto y$, $z\mapsto z-qu,t\mapsto t$ on two copies of $U$ and as the identity on the other summands of $L$. It is easily checked that these are isometries and that
\begin{align*} \varphi\cdot f(s,r)&=f(s+p,r), \\
\psi\cdot f(s,r)&=f(s,r+q).
\end{align*}
This implies that the image of $f$ is invariant under the action of $\Gamma:=\langle\varphi,\psi\rangle$, and also that the induced action on $\R^2$ is as desired. 

To show that the image of $f$ is contained in $K\Omega$, let $f_1,f_2,f_3$ be the three ``components'' of $f$, as before, and let $(s,r)\in\R^2$. We compute
\begin{align*}
(f_2(s,r),f_2(s,r))&=(x-su+2y+a,x-su+2y+a) \\
&=4(x,y)+(a,a) \\
&=4+(e,e)\geq 3\frac12>0, \\
(f_3(s,r),f_3(s,r))&=(z-ru+2t+b,z-ru+2t+b) \\
&=4(z,t)+(b,b) \\
&=4+(e,e)\geq 3\frac12>0, \\
(f_2(s,r),f_3(s,r))&=(x-su+2y+a,z-ru+2t+b) \\
&=0
\end{align*}
and conclude that $[f_2(s,r)+if_3(s,r)]\in\Omega$. Also,
\begin{align*}
(f_1(s,r),f_1(s,r))&=(2u+v+sy+rt,2u+v+sy+rt)=(2u,v)+(v,2u)=4>0 , \\
(f_1(s,r),f_2(s,r))&=(2u+v+sy+rt,x-su+2y+a) \\
&=-s(u,v)+s(x,y)=-s+s=0, \\
(f_1(s,r),f_3(s,r))&=(2u+v+sy+rt,z-ru+2t+b)=-r(u,v)+r(z,t)=-r+r=0
\end{align*}
implies that $f(s,r)\in K\Omega$. 

It remains to check that $f(s,r)\in K\Omega^0$ for all $(s,r)\in\R^2$. 

\begin{prof}
Let $d\in L$ such that $(d,d)=-2$ and $(d,f_1(s,r))=(d,f_2(s,r))=(d,f_3(s,r))=0$ and as before write
\[ d=Au+Bv+Cx+Dy+Ez+Ft+d_1+d_2,
\]
and $(d_i,d_i)=-2n_i$ for $n_i\in\Z_{\geq0}$.
We need to find a contradiction. The relevant equations now become 
\begin{align}
AB+CD+EF&=n_1+n_2-1,\label{eq:t2standard:1} \\ 
2B+A+Cs+Er&=0,\label{eq:t2standard:2} \\ 
D-Bs+2C+(d_2,e)&=0, \label{eq:t2standard:3} \\ 
F-Br+2E+(d_1,e)&=0. \label{eq:t2standard:4}
\end{align}
\smallskip

\textbf{Case $B=0$:} The assumptions on $e$, together with (\ref{eq:t2standard:3}) and (\ref{eq:t2standard:4}), imply that $D+2C=F+2E=0$ and $d_1=d_2=0$, so that $n_1=n_2=0$. But then (\ref{eq:t2standard:1}) becomes
\[ 2C^2+2E^2=1,
\]
which is impossible.
\smallskip 

\textbf{Case $B\neq0$:} From (\ref{eq:t2standard:3}) and (\ref{eq:t2standard:4}) we get
\begin{align*}
s&=\frac{D+2C+(d_2,e)}{B}, & r&=\frac{F+2E+(d_1,e)}{B}.
\end{align*}
Substituting this into (\ref{eq:t2standard:2}) gives
\[ AB+CD+EF=-2B^2-2C^2-2E^2-C(d_2,e)-E(d_1,e),
\]
and combining this with (\ref{eq:t2standard:1}) we obtain
\[ 2B^2+2C^2+2E^2+C(d_2,e)+E(d_1,e)=1-n_1-n_2.
\]
The assumptions on $e$ imply that $Cd_2+Ed_1=0$, so that this becomes
\[ 2B^2+2C^2+2E^2=1-n_1-n_2.
\]
This is impossible under the assumption $B\neq0$, since $n_i\in\Z_{\geq0}$.
\end{prof}

\subsection{The PMSCTs with leaf space a torus of type (\ref{torusII})}\label{example:t2weird}

In this example we will construct a PMSCT whose leaf space is a torus with an induced integral affine structure of type (\ref{torusII}), namely one induced by the action of $\Z^2$ on $\R^2$ generated by $(x,y)\mapsto(x+p,y)$ and $(x,y)\mapsto(x+ny,y+q)$, where $n,p,q\in\Z_{\geq1}$. 

Consider the map $f:\R^2\to L_{\R}\times\P(L_{\C})$ defined by
\[ (s,r)\mapsto (2u+v+sy+rt,[qx+(nr^2-qs)u-nrz+2qy+a+i(z-ru+2q^2t+2nqry+b)]),
\]
the map $\varphi:L\to L$ defined as before and the map $\psi:L\to L$ defined by $u\mapsto u,v\mapsto v+qt$, $x\mapsto x-nz+qnu,y\mapsto y$, $z\mapsto z-qu,t\mapsto t+ny$ on the copies of $U$ and the identity on the other summands of $L$. It is easily checked that these are isometries and that
\begin{align*} \varphi\cdot f(s,r)&=f(s+p,r), \\
\psi\cdot f(s,r)&=f(s+nr,r+q).
\end{align*}
This implies that the image of $f$ is invariant under the action of $\Gamma:=\langle\varphi,\psi\rangle$, and also that the induced action on $\R^2$ is the desired one. To show that the image of $f$ is contained in $K\Omega$, denote once more by $f_1,f_2,f_3$ the ``components'' of $f$, and let $(s,r)\in\R^2$. Since 
\begin{align*}
(f_2(s,r),f_2(s,r))&=(qx+(nr^2-qs)u-nrz+2qy+a,qx+(nr^2-qs)u-nrz \\
& \hspace{22.8em} +2qy+a) \\
&=4q^2(x,y)+(a,a) \\
&=4q^2+(e,e)\geq 3\frac12>0, \\
(f_3(s,r),f_3(s,r))&=(z-ru+2q^2t+2nqry+b,z-ru+2q^2t+2nqry+b) \\
&=4q^2(z,t)+(b,b) \\
&=4q^2+(e,e)\geq 3\frac12>0, \\
(f_2(s,r),f_3(s,r))&=(qx+(nr^2-qs)u-nrz+2qy+a,z-ru+2q^2t+2nqry+b) \\
&=2nq^2r(x,y)-2nq^2r(z,t)=2nq^2r-2nq^2r=0.
\end{align*}
we get that $[f_2(s,r)+if_3(s,r)]\in\Omega$. The computations
\begin{align*}
(f_1(s,r),f_1(s,r))&=(2u+v+sy+rt,2u+v+sy+rt)=(2u,v)+(v,2u)=4>0 , \\
(f_1(s,r),f_2(s,r))&=(2u+v+sy+rt,qx+(nr^2-qs)u-nrz+2qy+a) \\
&=(nr^2-qs)(u,v)+qs(x,y)-nr^2(z,t)=nr^2-qs+qs-nr^2=0, \\
(f_1(s,r),f_3(s,r))&=(2u+v+sy+rt,z-ru+2q^2t+2nqry+b) \\
&=-r(u,v)+r(z,t)=-r+r=0
\end{align*}
show that $f(s,r)\in K\Omega$. 

It remains to show that $f(s,r)\in K\Omega^0$ for all $(s,r)\in\R^2$. 
\begin{prof} Let $d\in L$ such that $(d,d)=-2$ and $(d,f_1(s))=(d,f_2(s))=(d,f_3(s))=0$. Like before we write
\[ d=Au+Bv+Cx+Dy+Ez+Ft+d_1+d_2,
\]
and we set $(d_i,d_i)=-2n_i$ with $n_i\in\Z_{\geq0}$. The goal is to find a contradiction. The main equations are now
\begin{align}
AB+CD+EF&=n_1+n_2-1,\label{eq:t2all:1} \\ 
2B+A+Cs+Er&=0,\label{eq:t2all:2} \\ 
Dq+B(nr^2-qs)-Fnr+2Cq+(d_2,e)&=0, \label{eq:t2all:3} \\ 
F-Br+2Eq^2+2Cnqr+(d_1,e)&=0. \label{eq:t2all:4}
\end{align}

\textbf{Case $B-2Cnq=0$:} Equation (\ref{eq:t2all:4}) tells us that $d_1=0$ and $F+2Eq^2=0$. 
\smallskip

\indent\indent\textbf{Subcase $C=0$:} This implies that $B=0$, so that (\ref{eq:t2all:1}) becomes
\[ 2E^2q^2=1-n_2.
\]
This is only possible if $E=0$ and $n_2=1$, but then also $F=0$ and (\ref{eq:t2all:3}) becomes
\[ Dq+(d_2,e)=0,
\]
which would imply that $d_2=0$, contradicting $n_2=1$.
\smallskip

\indent\indent\textbf{Subcase $C\neq0$:} Equation (\ref{eq:t2all:2}) tells us that
\[ s=-\frac{2B+A+Er}{C},
\]
and with (\ref{eq:t2all:3}) we obtain
\[ 2Cn^2qr^2-2Fnr+2Anq^2+C(8n^2q^3+2q)+Dq+(d_2,e)=0.
\]
Since $C,n,q\neq0$ and $r\in\R$, we must have that
\[ F^2\geq 2Cq\big[2Anq^2+C(8n^2q^3+2q)+Dq+(d_2,e)\big].
\]
But $F=-2Eq^2$ and $B=2Cnq$, so combining this with (\ref{eq:t2all:1}) yields
\[ q(1-n_2)\geq C^2(8n^2q^3+2q)+C(d_2,e).
\]
Since $C\neq0$, this is certainly impossible when $C$ and $(d_2,e)$ have the same parity. So let us assume that they have opposite parity, so that the equation becomes
\begin{equation}
q(1-n_2)\geq C^2(8n^2q^3+2q)-|C|\cdot|(d_2,e)|. \label{eq:ineq:all}
\end{equation}
Now both $d_2$ and $e$ lie in the same copy of $-E_8$, and since $(\cdot,\cdot)$ is negative definite on $-E_8$ we can use the Cauchy-Schwarz inequality to obtain
\[ |(d_2,e)|\leq\sqrt{|(d_2,d_2)|\cdot|(e,e)|}=\sqrt{2\cdot|(e,e)|n_2}\leq\sqrt{n_2},
\]
using that we chose $e$ such that $|(e,e)|\leq\frac12$. Now, in order for (\ref{eq:ineq:all}) to hold we certainly must have
\[ C^2(8n^2q^3+2q)-\sqrt{n_2}\cdot|C|+qn_2-q\leq0
\]
and it is easily seen that this is not possible for $0\neq C\in\Z$.
\smallskip

\textbf{Case $B-2Cnq\neq0$:} We immediately distinguish two cases: $B=0$ and $B\neq0$.
\smallskip

\indent\indent\textbf{Subcase $B=0$:} We claim that $F\neq0$. Indeed, if we had $F=0$, (\ref{eq:t2all:3}) would become
\[ Dq+2Cq+(d_2,e)=0,
\]
meaning that $d_2=0$, so $n_2=0$, and $D+2C=0$. But then (\ref{eq:t2all:1}) becomes
\[ 2C^2=1-n_1,
\]
which can only hold if $C=0$ and $n_1=1$. But then (\ref{eq:t2all:4}) becomes
\[ 2Eq^2+(d_1,e)=0,
\]
which implies $d_1=0$, contradicting $n_1=1$. So we see indeed that $F\neq0$. But then (\ref{eq:t2all:3}) and (\ref{eq:t2all:4}) yield
\[ r=-\frac{F+2Eq^2+(d_1,e)}{2Cnq}=\frac{Dq+2Cq+(d_2,e)}{Fn}.
\]
This becomes
\[ 2CDnq^2+4C^2nq^2+2Cnq(d_2,e)+F^2n+2EFnq^2+Fn(d_1,e)=0,
\]
and the assumptions on $e$ imply that $2Cqd_2+Fd_1=0$ and
\[ 2CDnq^2+4C^2nq^2+F^2n+2EFnq^2=0.
\]
Since $B=0$, combining this with (\ref{eq:t2all:1}) we obtain
\[ 4C^2nq^2+F^2n=2nq^2(1-n_1-n_2).
\]
Both $C$ and $F$ are nonzero, meaning that this is impossible.
\smallskip

\indent\indent\textbf{Subcase $B\neq0$:} We can write
\begin{align*}
r&=\frac{F+2Eq^2+(d_1,e)}{B-2Cnq}, & s&=\frac{Dq+Bnr^2-Fnr+2Cq+(d_2,e)}{Bq}.
\end{align*}
This yields
\begin{align*} s&=\frac{(B-2Cnq)^2(2Cq+Dq+(d_2,e))+Bn(F+2Eq^2+(d_1,e))^2}{Bq(B-2Cnq)^2}\\
&\ \ \ \ -\frac{Fn(B-2Cnq)(F+2Eq^2+(d_1,e))}{Bq(B-2Cnq)^2}
\end{align*}
and substituting this into (\ref{eq:t2all:2}) and using the assumptions on $e$ (actually, finally using them to their full potential), this reduces to
\begin{align*}
0&=2B^2q(B-2Cnq)^2+ABq(B-2Cnq)^2 \\
&\,\,\,\,\,\,\, +C\Big((B-2Cnq)^2(2Cq+Dq)+Bn(F+2Eq^2)^2-Fn(B-2Cnq)(F+2Eq^2)\Big) \\
&\,\,\,\,\,\,\, +BEq(B-2Cnq)(F+2Eq^2).
\end{align*}
Some rewriting turns this into
\begin{align*}
0&=q(B-2Cnq)^2\big(2B^2+2C^2+AB+CD\big) \\
&\,\,\,\,\,\,\, +BCn(F+2Eq^2)^2-CFn(B-2Cnq)(F+2Eq^2)+BEq(B-2Cnq)(F+2Eq^2),
\end{align*}
and some easy computations show that the second line is equal to
\[ EFq(B-2Cnq)^2+2q(BEq+CFn)^2,
\]
so that altogether we obtain
\[ q(B-2Cnq)^2\big(2B^2+2C^2+AB+CD+EF\big)+2q(BEq+CFn)^2=0.
\]
Combining this with (\ref{eq:t2all:1}) we get 
\[ 2\Big((B-2Cnq)^2(B^2+C^2)+(BEq+CFn)^2\Big)=(B-2Cnq)^2(1-n_1-n_2).
\]
But this is impossible, since $B-2Cnq\neq0$, $B\neq0$ and $n_1,n_2\geq0$, giving us the desired contradiction.
\end{prof}

\bibliographystyle{alpha}

\bibliography{references}

\begin{thebibliography}{{Mar}13}

\bibitem[ACH95]{AlcaldeCuesta1994IntegrationSD}
F.~Alcalde-Cuesta and G.~Hector.
\newblock Int\'egration symplectique des vari\'et\'es de poisson
  r\'eguli\`eres.
\newblock {\em Israel Journal of Math}, 90(1-3):125--165, 1995.

\bibitem[BPV84]{barth1984compact}
W.~Barth, C.~Peters, and A.~Ven.
\newblock {\em Compact complex surfaces}.
\newblock Ergebnisse der Mathematik und ihrer Grenzgebiete. Springer, 1984.

\bibitem[CF04]{intpoiss}
M.~Crainic and R.~L. Fernandes.
\newblock Integrability of {P}oisson brackets.
\newblock {\em J. Differential Geom.}, 66(1):71--137, 01 2004.

\bibitem[CF05]{rigidity}
M.~Crainic and R.~L. Fernandes.
\newblock Rigidity and flexibility in {P}oisson geometry.
\newblock In {\em Travaux math\'{e}matiques. {F}asc. {XVI}}, volume~16 of {\em
  Trav. Math.}, pages 53--68. Univ. Luxemb., Luxembourg, 2005.

\bibitem[CFM19]{pmct2}
M.~Crainic, R.~L. Fernandes, and D.~{Mart{\'i}nez Torres}.
\newblock Regular {P}oisson manifolds of compact types ({PMCT} 2).
\newblock {\em Asterisque}, 413:1--166, January 2019.

\bibitem[CFM21]{lectures}
M.~Crainic, R.~L. Fernandes, and I.~M\u{a}rcu\c{t}.
\newblock {\em Lectures on {P}oisson geometry}, volume 217 of {\em Graduate
  Studies in Mathematics}.
\newblock American Mathematical Society, Providence, RI, [2021] \copyright
  2021.

\bibitem[CFT19]{pmct1}
M.~Crainic, R.~L. Fernandes, and D.~Martínez Torres.
\newblock Poisson manifolds of compact types ({PMCT} 1).
\newblock {\em Journal für die reine und angewandte Mathematik}, 2019(756):101
  -- 149, 2019.

\bibitem[DH82]{dh}
J.~J. Duistermaat and G.~J. Heckman.
\newblock On the variation in the cohomology of the symplectic form of the
  reduced phase space.
\newblock {\em Inventiones mathematicae}, 69(2):259--268, 1982.

\bibitem[Gol22]{Goldman22}
W.~M. Goldman.
\newblock {\em Geometric structures on manifolds}, volume 227 of {\em Graduate
  Studies in Mathematics}.
\newblock American Mathematical Society, Providence, RI, [2022] \copyright
  2022.

\bibitem[Kot06]{kotschick2006free}
D.~Kotschick.
\newblock Free circle actions with contractible orbits on symplectic manifolds.
\newblock {\em Mathematische Zeitschrift}, 252(1):19--25, 2006.

\bibitem[{Mar}13]{2013arXiv1312.7267M}
D.~{Mart{\'{\i}}nez Torres}.
\newblock {A Poisson manifold of strong compact type}.
\newblock {\em Indagationes Mathematicae}, 25(5):1154--1159, December 2013.

\bibitem[Mis96]{MISHACHEV1996301}
K.N. Mishachev.
\newblock The classification of {L}agrangian bundles over surfaces.
\newblock {\em Differential Geometry and its Applications}, 6(4):301 -- 320,
  1996.

\bibitem[Mor53]{pjm/1103051334}
L.~J. Mordell.
\newblock {On the linear independence of algebraic numbers.}
\newblock {\em Pacific Journal of Mathematics}, 3(3):625 -- 630, 1953.

\bibitem[Sep10]{Kleinintaff}
D.~Sepe.
\newblock Classification of {L}agrangian fibrations over a {K}lein bottle.
\newblock {\em Geom. Dedicata}, 149:347--362, 2010.

\bibitem[Sep13]{strongintaff}
D.~Sepe.
\newblock Universal {L}agrangian bundles.
\newblock {\em Geometriae Dedicata}, 165, August 2013.

\bibitem[Siu83]{Siu1983}
Y.-T. Siu.
\newblock Every {K}3 surface is {K}{\"a}hler.
\newblock {\em Inventiones mathematicae}, 73(1):139--150, Feb 1983.

\end{thebibliography}

\end{document}